# On likelihood ratio tests

## Erich L. Lehmann[1]

*University of California at Berkeley*

**Abstract:** Likelihood ratio tests are intuitively appealing. Nevertheless, a number of examples are known in which they perform very poorly. The present paper discusses a large class of situations in which this is the case, and analyzes just how intuition misleads us; it also presents an alternative approach which in these situations is optimal.

## 1. The popularity of likelihood ratio tests

Faced with a new testing problem, the most common approach is the likelihood ratio (LR) test. Introduced by Neyman and Pearson in 1928, it compares the maximum likelihood under the alternatives with that under the hypothesis. It owes its popularity to a number of facts.

(i) It is intuitively appealing. The likelihood of $\theta$,

$$L_x(\theta) = p_\theta(x)$$

i.e. the probability density (or probability) of $x$ considered as a function of $\theta$, is widely considered a (relative) measure of support that the observation $x$ gives to the parameter $\theta$. (See for example Royall [8]). Then the likelihood ratio

$$\sup_{alt}[p_\theta(x)] \big/ \sup_{hyp}[p_\theta(x)] \tag{1.1}$$

compares the best explanation the data provide for the alternatives with the best explanations for the hypothesis. This seems quite persuasive.

(ii) In many standard problems, the LR test agrees with tests obtained from other principles (for example it is UMP unbiased or UMP invariant). Generally it seems to lead to satisfactory tests. However, counter-examples are also known in which the test is quite unsatisfactory; see for example Perlman and Wu [7] and Menéndez, Rueda, and Salvador [6].

(iii) The LR test, under suitable conditons, has good asymptotic properties.

None of these three reasons are convincing.

(iii) tells us little about small samples.

(i) has no strong logical grounding.

(ii) is the most persuasive, but in these standard problems (in which there typically exist a complete set of sufficient statistics) all principles typically lead to tests that are the same or differ only by little.

[1]Department of Statistics, 367 Evans Hall, University of California, Berkeley, CA 94720-3860, e-mail: shaffer@stat.berkeley.edu

In view of lacking theoretical support and many counterexamples, it would be good to investigate LR tests systematically for small samples, a suggestion also made by Perlman and Wu [7]. The present paper attempts a first small step in this endeavor.

## 2. The case of two alternatives

The simplest testing situation is that of testing a simple hypothesis against a simple alternative. Here the Neyman-Pearson Lemma completely vindicates the LR-test, which always provides the most powerful test. Note however that in this case no maximization is involved in either the numerator or denominator of (1.1), and as we shall see, it is just these maximizations that are questionable.

The next simple situation is that of a simple hypothesis and two alternatives, and this is the case we shall now consider.

Let $\underline{X}=(X_1, \ldots, X_n)$ where the $X$'s are iid. Without loss of generality suppose that under $H$ the $X$'s are uniformly distributed on $(0,1)$. Consider two alternatives $f, g$ on $(0,1)$. To simplify further, we shall assume that the alternatives are symmetric, i.e. that

$$p_1(\underline{x}) = f(x_1) \cdots f(x_n)$$
$$p_2(\underline{x}) = f(1-x_1) \cdots f(1-x_n). \tag{2.1}$$

Then it is natural to restrict attention to symmetric tests (that is the invariance principle) i.e. to rejection regions $R$ satisfying

$$(x_1, \ldots, x_n) \in R \quad \text{if and only if } (1-x_1, \ldots, 1-x_n) \in R. \tag{2.2}$$

The following result shows that under these assumptions there exists a uniformly most powerful (UMP) invariant test, i.e. a test that among all invariant tests maximizes the power against both $p_1$ and $p_2$.

**Theorem 2.1.** *For testing $H$ against the alternatives (2.1) there exists among all level $\alpha$ rejection regions $R$ satisfying (2.2) one that maximizes the power against both $p_1$ and $p_2$ and it rejects $H$ when*

$$\frac{1}{2}[p_1(\underline{x}) + p_2(\underline{x})] \quad \text{is sufficiently large.} \tag{2.3}$$

We shall call the test (2.3) the *average likelihood ratio test* and from now on shall refer to (1.1) as the *maximum likelihood ratio test*.

*Proof.* If $R$ satisfies (2.2), its power against $p_1$ and $p_2$ must be the same. Hence

$$\int_R p_1 \quad = \quad \int_R p_2 = \int_R \frac{1}{2}(p_1 + p_2). \tag{2.4}$$

$\square$

By the Neyman–Pearson Lemma, the most powerful test of $H$ against $\frac{1}{2}[p_1 + p_2]$ rejects when (2.3) holds.

**Corollary 2.1.** *Under the assumptions of Theorem 2.1, the average LR test has power greater than or equal to that of the maximum likelihood ratio test against both $p_1$ and $p_2$.*

*Proof.* The maximum LR test rejects when

$$\max(p_1(\underline{x}),\ p_2(\underline{x})) \quad \text{is sufficiently large.} \tag{2.5}$$

Since this test satisfies (2.2), the result follows. □

The Corollary leaves open the possibility that the average and maximum LR tests have the same power; in particular they may coincide. To explore this possibility consider the case $n = 1$ and suppose that $f$ is increasing. Then the likelihood ratio will be

$$f(x) \quad \text{if } x > \tfrac{1}{2} \quad \text{and } f(1-x) \text{ if } x < \tfrac{1}{2}. \tag{2.6}$$

The maximum LR test will therefore reject when

$$\left| x - \frac{1}{2} \right| \quad \text{is sufficiently large} \tag{2.7}$$

i.e. when $x$ is close to either 0 or 1.

It turns out that the average LR test will depend on the shape of $f$ and we shall consider two cases: (a) $f$ is convex; (b) $f$ is concave.

**Theorem 2.2.** *Under the assumptions of Theorem 2.1 and with $n = 1$,*

(i) (a) *if $f$ is convex, the average LR test rejects when (2.7) holds;*
  (b) *if $f$ is concave, the average LR test rejects when*

$$\left| x - \frac{1}{2} \right| \quad \text{is sufficiently small.} \tag{2.8}$$

(ii) (a) *if $f$ is convex, the maximum LR test coincides with the average LR test, and hence is UMP among all tests satisfying (2.2) for n=1.*
  (b) *if $f$ is concave, the maximum LR test uniformly* minimizes *the power among all tests satisfying (2.2) for n=1, and therefore has power $< \alpha$.*

*Proof.* This is an immediate consequence of the fact that if $x < x' < y' < y$ then

$$\frac{f(x)+f(y)}{2} \quad \text{is} \quad \genfrac{}{}{0pt}{}{>}{<} \quad \frac{f(x')+f(y')}{2} \quad \text{if } f \text{ is} \quad \genfrac{}{}{0pt}{}{\text{convex.}}{\text{concave.}} \tag{2.9}$$

It is clear from the argument that the superiority of the average over the likelihood ratio test in the concave case will hold even if $p_1$ and $p_2$ are not exactly symmetric. Furthermore it also holds if the two alternatives $p_1$ and $p_2$ are replaced by the family $\theta p_1 + (1-\theta)p_2$, $0 \le \theta \le 1$. □

## 3. A finite number of alternatives

The comparison of maximum and average likelihood ratio tests discussed in Section 2 for the case of two alternatives obtains much more generally. In the present section we shall sketch the corresponding result for the case of a simple hypothesis against a finite number of alternatives which exhibit a symmetry generalizing (2.1).

Suppose the densities of the simple hypothesis and the $s$ alternatives are denoted by $p_0, p_1, \ldots, p_s$ and that there exists a group $G$ of transformations of the sample which leaves invariant both $p_0$ and the set $\{p_1, \ldots, p_s\}$ (i.e. each transformation results in a permutation of $p_1, \ldots, p_s$). Let $\bar{G}$ denote the set of these permutations and suppose that it is transitive over the set $\{p_1, \ldots, p_s\}$ i.e. that given any $i$ and $j$ there exists a transformation in $\bar{G}$ taking $p_i$ into $p_j$. A rejection region $R$ is said to be invariant under $\bar{G}$ if

$$x \in R \quad \text{if and only if} \quad g(x) \in R \text{ for all } g \text{ in } G. \qquad (3.1)$$

**Theorem 3.1.** *Under these assumptions there exists a uniformly most powerful invariant test and it rejects when*

$$\frac{\sum_{i=1}^s p_i(x)/s}{p_0(x)} \quad \text{is sufficiently large.} \qquad (3.2)$$

In generalization of the terminology of Theorem 2.1 we shall call (3.2) the average likelihood ratio test. The proof of Theorem 3.1 exactly parallels that of Theorem 2.1.

The Theorem extends to the case where $G$ is a compact group. The average in the numerator of (3.2) is then replaced by the integral with respect to the (unique) invariant probability measure over $\bar{G}$. For details see Eaton ([3], Chapter 4). A further extension is to the case where not only the alternatives but also the hypothesis is composite.

To illustrate Theorem 3.1, let us extend the case considered in Section 2. Let $(X, Y)$ have a bivariate distribution over the unit square which is uniform under $H$. Let $f$ be a density for $(X, Y)$ which is strictly increasing in both variables and consider the four alternatives

$$p_1 = f(x, y), \, p_2 = f(1-x, y), \, p_3 = f(x, 1-y), \, p_4 = f(1-x, 1-y).$$

The group $G$ consists of the four transformations

$g_1(x, y) = (x, y), \, g_2(x, y) = (1-x, y), \, g_3(x, y) = (x, 1-y),$
and $g_4(x, y) = (1-x, 1-y).$

They induce in the space of $(p_1, \ldots, p_4)$ the transformations:

$\bar{g}_1 = $ the identity
$\bar{g}_2: p_1 \rightarrow p_2, \, p_2 \rightarrow p_1, \, p_3 \rightarrow p_4, \, p_4 \rightarrow p_3$
$\bar{g}_3: p_1 \rightarrow p_3, \, p_3 \rightarrow p_1, \, p_2 \rightarrow p_4, \, p_4 \rightarrow p_2$
$\bar{g}_4: p_1 \rightarrow p_4, \, p_4 \rightarrow p_1, \, p_2 \rightarrow p_3, \, p_3 \rightarrow p_2.$

This is clearly transitive, so that Theorem 3.1 applies. The uniformly most powerful invariant test, which rejects when

$$\sum_{i=1}^4 p_i(x, y) \quad \text{is large}$$

is therefore uniformly at least as powerful as the maximum likelihood ratio test which rejects when

$$\max\left[p_1\left(x, y\right), p_2\left(x, y\right), p_3\left(x, y\right), p_4\left(x, y\right)\right]$$

is large.

## 4. Location-scale families

In the present section we shall consider some more classical problems in which the symmetries are represented by infinite groups which are not compact. As a simple example let the hypothesis $H$ and the alternatives $K$ be specified respectively by

$$H: f(x_1 - \theta, \ldots, x_n - \theta) \text{ and } K: g(x_1 - \theta, \ldots, x_n - \theta) \quad (4.1)$$

where $f$ and $g$ are given densities and $\theta$ is an unknown location parameter. We might for example want to test a normal distribution with unknown mean against a logistic or Cauchy distribution with unknown center.

The symmetry in this problem is characterized by the invariance of $H$ and $K$ under the transformations

$$X_i' = X_i + c \ (i = 1, \ldots, n). \quad (4.2)$$

It can be shown that there exists a uniformly most powerful invariant test which rejects $H$ when

$$\frac{\int_{-\infty}^{\infty} g(x_1 - \theta, \ldots, x_n - \theta) d\theta}{\int_{-\infty}^{\infty} f(x_1 - \theta, \ldots, x_n - \theta) d\theta} \quad \text{is large.} \quad (4.3)$$

The method of proof used for Theorem 2.1 and which also works for Theorem 3.1 no longer works in the present case since the numerator (and denominator) no longer are averages. For the same reason the term average likelihood ratio is no longer appropriate and is replaced by integrated likelihood. However an easy alternative proof is given for example in Lehmann ([5], Section 6.3).

In contrast to (4.2), the maximum likelihood ratio test rejects when

$$\frac{g(x_1 - \hat{\theta}_1, \ldots, x_n - \hat{\theta}_1)}{f(x_1 - \hat{\theta}_0, \ldots, x_n - \hat{\theta}_0)} \quad \text{is large,} \quad (4.4)$$

where $\hat{\theta}_1$ and $\hat{\theta}_0$ are the maximum likelihood estimators of $\theta$ under $g$ and $f$ respectively. Since (4.4) is also invariant under the transformations (4.2), it follows that the test (4.3) is uniformly at least as powerful as (4.4), and in fact more powerful unless the two tests coincide which will happen only in special cases.

The situation is quite similar for scale instead of location families. The problem (4.1) is now replaced by

$$H: \frac{1}{\tau^n} f\left(\frac{x_1}{\tau}, \ldots, \frac{x_n}{\tau}\right) \quad \text{and} \quad K: \frac{1}{\tau^n} g\left(\frac{x_1}{\tau}, \ldots, \frac{x_n}{\tau}\right) \quad (4.5)$$

where either the $x$'s are all positive or $f$ and $g$ symmetric about 0 in each variable.

This problem remains invariant under the transformations

$$X_i' = cX_i, \ c > 0. \quad (4.6)$$

It can be shown that a uniformly most powerful invariant test exists and rejects $H$ when

$$\frac{\int_0^{\infty} \nu^{n-1} g(\nu x_1, \ldots, \nu x_n) d\nu}{\int_0^{\infty} \nu^{n-1} f(\nu x_1, \ldots, \nu x_n) d\nu} \quad \text{is large.} \quad (4.7)$$

On the other hand, the maximum likelihood ratio test rejects when

$$\frac{g(\frac{x_1}{\hat{\tau}_1}, \dots, \frac{x_n}{\hat{\tau}_n})/\hat{\tau}_1^n}{f(\frac{x_1}{\hat{\tau}_0}, \dots, \frac{X_n}{\hat{\tau}_0.})/\hat{\tau}_0^n} \quad \text{is large} \tag{4.8}$$

where $\hat{\tau}_1$ and $\hat{\tau}_0$ are the maximum likelihood estimators of $\tau$ under $g$ and $f$ respectively. Since it is invariant under the transformations (4.6), the test (4.8) is less powerful than (4.7) unless they coincide.

As in (4.3), the test (4.7) involves an integrated likelihood, but while in (4.3) the parameter $\theta$ was integrated with respect to Lebesgue measure, the nuisance parameter in (4.6) is integrated with respect to $\nu^{n-1}d\nu$. A crucial feature which all the examples of Sections 2–4 have in common is that the group of transformations that leave $H$ and $K$ invariant is transitive i.e. that there exists a transformation which for any two members of $H$ (or of $K$) takes one into the other. A general theory of this case is given in Eaton ([3], Sections 6.7 and 6.4).

Elimination of nuisance parameters through integrated likelihood is recommended very generally by Berger, Liseo and Wolpert [1]. For the case that invariance considerations do not apply, they propose integration with respect to non-informative priors over the nuisance parameters. (For a review of such prior distributions, see Kass and Wasserman [4]).

## 5. The failure of intuition

The examples of the previous sections show that the intuitive appeal of maximum likelihood ratio tests can be misleading. (For related findings see Berger and Wolpert ([2], pp. 125–135)). To understand just how intuition can fail, consider a family of densities $p_\theta$ and the hypothesis $H: \theta = 0$. The Neyman–Pearson lemma tells us that when testing $p_0$ against a specific $p_\theta$, we should reject $y$ in preference to $x$ when

$$\frac{p_\theta(x)}{p_0(x)} \quad < \quad \frac{p_\theta(y)}{p_0(y)} \tag{5.1}$$

the best test therefore rejects for large values of $p_\theta(x)/p_0(x)$, i.e. is the maximum likelihood ratio test.

However, when more than one value of $\theta$ is possible, consideration of only large values of $p_\theta(x)/p_0(x)$ (as is done by the maximum likelihood ratio test) may no longer be the right strategy. Values of $x$ for which the ratio $p_\theta(x)/p_0(x)$ is small now also become important; they may have to be included in the rejection region because $p_{\theta'}(x)/p_0(x)$ is large for some other value $\theta'$.

This is clearly seen in the situation of Theorem 2 with $f$ increasing and $g$ decreasing, as illustrated in Fig. 1.

For the values of $x$ for which $f$ is large, $g$ is small, and vice versa. The behavior of the test therefore depends crucially on values of $x$ for which $f(x)$ or $g(x)$ is small, a fact that is completely ignored by the maximum likelihood ratio test.

Note however that this same phenomenon does not arise when all the alternative densities $f, g, \dots$ are increasing. When $n = 1$, there then exists a uniformly most powerful test and it is the maximum likelihood ratio test. This is no longer true when $n > 1$, but even then all reasonable tests, including the maximum likelihood ratio test, will reject the hypothesis in a region where all the observations are large.

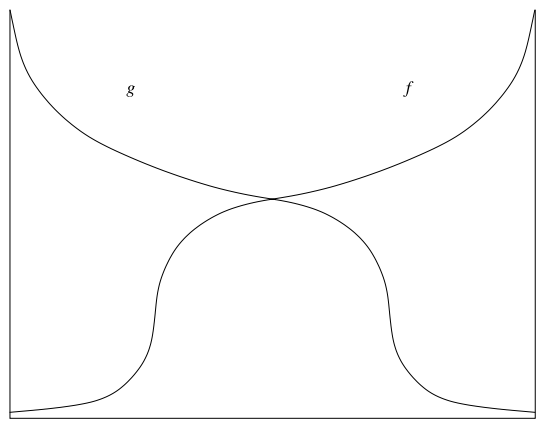

Fɪɢ 1.

## 6. Conclusions

For the reasons indicated in Section 1, maximum likelihood ratio tests are so widely accepted that they almost automatically are taken as solutions to new testing problems. In many situations they turn out to be very satisfactory, but gradually a collection of examples has been building up and is augmented by those of the present paper, in which this is not the case.

In particular when the problem remains invariant under a transitive group of transformations, a different principle (likelihood averaged or integrated with respect to an invariant measure) provides a test which is uniformly at least as good as the maximum likelihood ratio test and is better unless the two coincide. From the argument in Section 2 it is seen that this superiority is not restricted to invariant situations but persists in many other cases. A similar conclusion was reached from another point of view by Berger, Liseo and Wolpert [1].

The integrated likelihood approach without invariance has the disadvantage of not being uniquely defined; it requires the choice of a measure with respect to which to integrate. Typically it will also lead to more complicated test statistics. Nevertheless: In view of the superiority of integrated over maximum likelihood for large classes of problems, and the considerable unreliability of maximum likelihood ratio tests, further comparative studies of the two approaches would seem highly desirable.